\documentclass[12pt, a4paper]{amsart}
\usepackage{amsmath}
\usepackage{geometry,amsthm,graphics,tabularx,amssymb,shapepar}
\usepackage{amscd}
\usepackage{bm}
\usepackage[all]{xypic}
\newcommand{\fv}{\mathfrak v}
\newcommand{\ff}{\mathfrak f}
\newcommand{\fb}{\mathfrak b}

\newcommand{\fg}{\mathfrak g}
\newcommand{\fh}{\mathfrak h}
\newcommand{\dg}{\dot{\mathfrak g}}
\newcommand{\dfh}{\dot{\mathfrak h}}

\newcommand{\al}{\alpha}

\newcommand{\wt}{\widetilde}
\newcommand{\wh}{\widehat}

\newcommand{\dal}{\dot{\alpha}}

\newcommand{\ot}{\otimes}

\newcommand{\CH}{\mathcal{H}}

\newcommand{\C}{\mathbb{C}}

\newcommand{\N}{\mathbb N}
\newcommand{\Z}{\mathbb Z}
\newcommand{\rd}{\mathrm{d}}
\newcommand{\rk}{\mathrm{k}}

\newcommand{\la}{\langle}
\newcommand{\ra}{\rangle}

\newcommand{\ba}{\begin {eqnarray}}
\newcommand{\ea}{\end {eqnarray}}
\newcommand{\baa}{\begin {eqnarray*}}
\newcommand{\eaa}{\end {eqnarray*}}
\newcommand{\be}{\begin {equation}}
\newcommand{\ee}{\end {equation}}
\newcommand{\bee}{\begin {equation*}}
\newcommand{\eee}{\end {equation*}}

\newcommand{\U}{\mathcal{U}}

\newcommand{\te}[1]{\textnormal{{#1}}}

\theoremstyle{Theorem}

\theoremstyle{Theorem}

\newtheorem{thm}{Theorem}[section]

\newtheorem{lemt}[thm]{Lemma}
\newtheorem{prpt}[thm]{Proposition}

\theoremstyle{Theorem}

\theoremstyle{Theorem}

\theoremstyle{Plain}

\theoremstyle{Definition}

\def\({\left(}

\def\){\right)}

\def \<{{\langle}}
\def \>{{\rangle}}

\allowdisplaybreaks

\numberwithin{equation}{section}
\title[Repn of toroidal EALA]{Integrable representations for toroidal extended affine Lie algebras}

\author{Fulin Chen$^1$}
\thanks{{1. Partially supported by China NSF grants (No.11501478)}}
\address{Department of Mathematics, Xiamen University,
 Xiamen, China 361005} \email{chenf@xmu.edu.cn }
\author{Zhiqiang Li}
\address{Department of Mathematics, Xiamen University,
 Xiamen, China 361005} \email{lzq06031212@sina.com}
 \author{Shaobin Tan$^2$}
 \thanks{{2. Partially supported by China NSF grants (Nos.11471268, 11531004)}}
 \address{Department of Mathematics, Xiamen University,
 Xiamen, China 361005} \email{tans@xmu.edu.cn}


\begin{document}

\begin{abstract}
Let $\fg$ be any untwisted affine Kac-Moody algebra, $\mu$ any fixed complex number, and $\wt\fg(\mu)$ the corresponding toroidal extended affine Lie algebra of nullity two. For any $k$-tuple $\bm{\lambda}=({\lambda}_1, \cdots, {\lambda}_k)$ of weights of $\fg$,
and $k$-tuple $\bm{a}=(a_1,\cdots, a_k)$ of distinct non-zero complex numbers, we construct a class of modules
$\wt V(\bm{\lambda},\bm{a})$ for the extended affine Lie algebra $\wt\fg(\mu)$. We prove that the $\wt\fg(\mu)$-module $\wt V(\bm{\lambda},\bm{a})$ is completely reducible.
We also prove that the $\wt\fg(\mu)$-module $\wt V(\bm{\lambda},\bm{a})$ is integrable when all weights $\lambda_i$ in $\bm{\lambda}$ are dominant integral. Thus, we obtain a new class of irreducible integrable weight modules for the toroidal extended affine Lie algebra $\wt\fg(\mu)$.
\end{abstract}
\maketitle

\section{Introduction}
In this paper we study representations for the nullity $2$ toroidal extended affine Lie algebras (EALAs for short). The notion of EALAs was first introduced by Hoegh-Krohn and Torresani in \cite{H-KT} under the name of quasi-simple Lie algebras.
By definition an EALA is a complex Lie algebra, together with a nonzero finite dimensional
ad-diagonalizable subalgebra and a non-degenerate invariant symmetric bilinear form, that satisfies  a list of natural axioms. The structure theory of EALAs has been intensively studied for the past twenty years (see \cite{AABGP,BGK,N} and the references therein). We recall that the rank of the group generated by all isotropic roots of an EALA is called the nullity of the Lie algebra.
Indeed, the nullity $0$ EALAs are nothing but the finite dimensional simple Lie algebras,
and the affine Kac-Moody algebras are precisely the nullity $1$  EALAs \cite{ABGP}.
EALAs of nullity 2 are closely related to the Lie algebras studied by Saito and
Slodowy on elliptic singularities (see \cite{S}). We know that the toroidal EALAs give a class of important EALAs that provide examples with arbitrary  nullities. In this paper, motivated by  Chari-Pressley's loop module construction \cite{CP} for affine Kac-Moody algebras and Billig's construction \cite{B} for toroidal EALAs, we construct a new class of irreducible integrable modules for the nullity $2$ toroidal extended affine Lie algebras, which indeed generalizes Billig's construction.

We first recall  the loop module construction for affine Kac-Moody algebras given  in \cite{CP}.
Let $\dg$ be a finite dimensional simple Lie algebra over $\C$ with a fixed Cartan subalgebra $\dfh$, and let $\dfh^*$ be the dual space of $\dfh$.
We denote by
\begin{align}\label{affalg}\fg=(\C[t_0,t_0^{-1}]\ot \dg)\oplus \C\rk_0 \oplus \C\rd_0\end{align}
the untwisted affinization of $\dg$.
Let $k$ be a fixed positive integer. For any pair
\[(\dot{\bm{\lambda}},\bm{a})\in (\dfh^*)^k\times (\C^\times)^k,\quad \dot{\bm{\lambda}}=(\dot{\lambda}_1,\cdots,\dot{\lambda}_k),\quad \bm{a}=(a_1,\cdots,a_k),\]
 we define  $V(\dot{\bm{\lambda}},\bm{a})$ to be the following loop vector space,
\begin{align}\label{vdlambdaa}
V_{\dg}(\dot{\lambda}_1)\ot\cdots\ot V_{\dg}(\dot{\lambda}_k)\ot \C[t_0,t_0^{-1}],\end{align}
which gives a module structure \cite{CP} for the affine Kac-Moody algebra $\fg$ with actions defined as follows
\begin{align}\label{aloopaction}
(t_0^n\ot x). v_1\ot \cdots \ot v_k\ot t_0^m=\sum_{i=1}^k a_i^n v_1\ot \cdots \ot (x.v_i)\ot \cdots \ot v_k\ot t_0^{m+n},
\end{align}
for $t_0^n\ot x$ $(n\in \Z, x\in \dg)$, and
 $\rk_0$ acts trivially, $\rd_0$ acts as differential operator $t_0\frac{d}{dt_0}$.
Here, the notation $V_{\dg}(\dot{\lambda}_i)$ stands for the irreducible highest weight $\dg$-module with highest weight $\dot{\lambda}_i$.
It is proved in \cite{CP} that, if the scalars $a_1,a_2,\cdots,a_k$ are distinct, then
the $\fg$-module $V(\dot{\bm{\lambda}},\bm{a})$ is completely reducible with finitely many irreducible components.
Moreover, if all $\dot\lambda_i\in \dfh^*$ are dominant weights, then every irreducible components of  $V(\dot{\bm{\lambda}},\bm{a})$  are integrable, and conversely, it is proved in \cite{C} that
every irreducible level $0$ integrable $\fg$-module raises in this way up to a shift action of $\rd_0$.

Let $\mathcal R=\C[t_0^{\pm 1}, t_1^{\pm 1}]$ be the ring of Laurent polynomials in commuting variables, and let $\mathcal S$ be the
subspace of divergence zero derivations on $\mathcal R$, which is also called set of skew-derivations \cite{BGK}.
Similar to the construction of untwisted affine Kac-Moody algebras, we can define a nullity $2$ toroidal EALA
\[\wt\fg=(\mathcal R\ot \dg)\oplus \mathcal K\oplus \mathcal S\]
 by adding the subspace of divergence zero derivations $\mathcal S$
to the
 universal central extension $(\mathcal R\ot \dg)\oplus \mathcal K$ of the iterated loop Lie algebra $\mathcal R\ot \dg$. More generally, we can add an abelian extension of $\mathcal S$ over $\mathcal K$ with a $2$-cocycle $\tau_\mu$, $\mu\in \C$, to get a more general toroidal EALAs $\wt\fg(\mu)$ (see Section 2 for details).

Comparing with their structure theory, the representation theory of EALAs with nullity greater than or equal to 2 are much less understood. So far,
the known  irreducible modules of toroidal EALAs were those constructed
in \cite{B} by applying the theory of vertex operator algebras.

Let $\fh=\dfh\oplus \C\rk_0\oplus \C\rd_0$ be the usual Cartan subalgebra of $\fg$, and $\fh^*$ the dual space of $\fh$. We let
\begin{align}\label{lambdaa}
(\bm{\lambda},\bm{a})\in (\fh^*)^k\times (\C^\times)^k,\quad \bm{\lambda}=(\lambda_1,\cdots,\lambda_k),\quad \bm{a}=(a_1,\cdots,a_k),\end{align}
be a pair such  that
\begin{align}\label{assumption}
\lambda_i(\rk_0)\ne 0\ \te{for all}\ i,\ \te{and}\ a_1,a_2,\cdots, a_k\ \te{are distinct}.\end{align}
  Motivated by Chari-Pressely's loop module construction for the affine Kac-Moody algebras, we construct a class of modules $\wt V(\bm{\lambda},\bm{a})$ for the nullity $2$ toroidal EALA $\wt\fg(\mu)$ with  fixed pair $(\bm{\lambda},\bm{a})$.
We remark that when $k=1$, it coincides with Billig's module construction \cite{B} for the nullity $2$ toroidal EALA $\wt\fg(\mu)$.
 Furthermore, we prove that the $\wt\fg(\mu)$-module $\wt V(\bm{\lambda},\bm{a})$ is completely reducible and has  finitely many
irreducible components. Thus we obtain in this way a new class of irreducible $\wt\fg(\mu)$-modules.

Let $\dot\fb=\C\rk_1\oplus \C\rd_1$ be an ablian Lie algebra equipped with a non-degenerate symmetric bilinear form $\<\cdot,\cdot\>$ determined by
$\<\rd_1,\rk_1\>=1$ and $\<\rd_1,\rd_1\>=0=\<\rk_1,\rk_1\>$.
We denote by
\begin{align}\label{viraff}
\bar{\ff}=\mathrm{Der}(\C[t_0,t_0^{-1}])\ltimes (\C[t_0,t_0^{-1}]\ot \dot{\ff}) \oplus \C\rk_0\oplus \C\rk_v\end{align}
 the so-called  Virasoro-affine Lie algebra associated to the reductive Lie algebra $\dot{\ff}=\dg\oplus \dot\fb$, and denote by
\begin{align}\label{barh}
\bar\ff_0=\dfh\oplus \C\rk_0\oplus \C\rd_0\oplus\C\rk_1\oplus \C\rd_1\oplus \C\rk_v=\fh\oplus \dot{\fb}\oplus \C\rk_v, \end{align}
a Cartan subalgebra of the Virasoro-affine Lie algebra $\bar\ff$.
For any $\lambda\in \fh^*$, we define a linear functional $\bar\lambda$ over $\bar{\ff}_0$ by letting
\begin{align}\label{barlambda}
\bar\lambda|_\fh=\lambda,\quad \bar\lambda(\rk_v)=24\mu \lambda(\rk_0),\quad \bar\lambda(\rk_1)=0=\bar\lambda(\rd_1).\end{align}
  Let $V_{\bar{\ff}}(\bar\lambda_i)$, $i=1,\cdots,k$, be the irreducible highest weight $\bar\ff$-module with highest weight $\bar{\lambda}_i$ (see Section 3). We now define a module $\wt V(\bm{\lambda},\bm{a})$ for the nullity $2$ toroidal EALA $\wt\fg(\mu)$ with a fixed pair $(\bm{\lambda},\bm{a})$. The underlying vector space of $\wt V(\bm{\lambda},\bm{a})$ is similar to loop module given by Chari-Pressely for the affine Kac-Moody algebras:
\begin{align}\label{loopmodspace}
\wt V(\bm{\lambda},\bm{a})=V_{\bar{\ff}}(\bar{\lambda}_1)\ot \cdots \ot V_{\bar{\ff}}(\bar{\lambda}_k) \ot \C[t_1,t_1^{-1}].
 \end{align}

 Since the Virasoro-affine Lie algebra $\bar\ff$ is not a Lie subalgebra of the nullity $2$ toroidal EALA $\wt\fg(\mu)$, the $\wt\fg(\mu)$-module action defined over the loop space $\wt V(\bm{\lambda},\bm{a})$ is very indirect and much more complicated then that defined by Chari-Pressely for the affine Kac-Moody algebra (see Section 4 for details). Indeed the action is constructed by a reformulation of Billig's technique (see Proposition \ref{refor}), and the major portions of the paper is devoted to the proof of the reducibility of the $\wt\fg(\mu)$-module $\wt V(\bm{\lambda},\bm{a})$.

We recall that one of the key steps in  Chari-Pressely's proof for the complete reducibility of $\fg$-module $V(\dot{\bm{\lambda}},\bm{a})$ is the Vandermonde determinant argument.
More precisely, by taking $n=1,2,\cdots,k$ and $m=l-n$ in \eqref{aloopaction} for some $l\in \Z$, there exists a system of linear equations whose coefficient matrix is a
Vandermonde matrix $(a_i^j)_{1\le i,j\le k}$, and therefore each vector
$v_1\ot \cdots\ot (x.v_i)\ot \cdots \ot v_k\ot t_0^l$ can be expressed as  a linear combination of the elements
$(t_0^n\ot x).v_1\ot \cdots\ot v_k\ot t^m$, $n,m\in \Z$.
However, this argument can not be directly applied to  our case. This is duo to that fact that the action of $t_1^n\ot x$, $x\in \dg$ on $V(\bm{\lambda},\bm{a})$ is defined by
\[(t_1^n\ot x). v_1\ot \cdots \ot v_i\ot \cdots\ot v_k\ot t_1^l=\sum_{i=1}^k a_i^n v_1\ot \cdots  \ot (x(0,n).v_i)\ot \cdots\ot v_k\ot t_1^{n+l}\]
for some operators $x(0,n)$ on $V_{\bar\ff}(\bar{\lambda}_i)$. While in the affine Kac-Moody algebra case (see \eqref{aloopaction}) the actions on the component $v_i$ is independent of $n$.

 In this paper the proof of the complete reducibility of the $\wt\fg(\mu)$-module $\wt V(\bm{\lambda},\bm{a})$  is based on a more subtle analysis of the operators such as $x(0,n)$. And a generalization of the classical Vandermonde determinant (see Lemma \ref{van}) will be applied to the refined argument.
 We also prove that, if all $\lambda_i$ are dominant integral weights, the $\wt\fg(\mu)$-module $\wt V(\bm{\lambda},\bm{a})$ is integrable and each weight spaces are finite dimensional. This result is clearly crucial in the study of the classification of irreducible integrable $\wt\fg(\mu)$-modules with finite dimensional weight spaces and non-zero central charges.

\vspace{0.2cm}
 Now we outline the structure of our paper. In Section 2 we briefly review some basic results we need for the nullity $2$ toroidal EALAs.
 In Section 3 we recall Billig's technique for the construction of modules for the EALAs, and give a reformulation of Billig's work.
 The main results of the paper, Theorem 4.5 and Proposition 4.7, are provided in Section 4. To verify the main results it is vital to prove a preparation result Theorem 4.2, and Section 5 is devoted to the proof of  Theorem \ref{irr}.

The sets of integers, non-negative integers, complex numbers and nonzero complex numbers will be denoted respectively by $\Z$, $\N$, $\C$ and $\C^\times$.

\section{Nullity $2$ toroidal EALAs}
In this section we recall the definition of  nullity $2$ toroidal EALA $\wt\fg(\mu)$ and review some basic results about the Lie algebras, which will be used later on.

We denote by $\Omega_\mathcal R^1$ the space of K\"{a}hler differentials on the ring $\mathcal R=\C[t_0^{\pm 1}, t_1^{\pm 1}]$.
As a left $\mathcal R$-module, $\Omega_\mathcal R^1$ has a basis $\rk_0, \rk_1$, where $\rk_i$, $i=0,1$, stands for the $1$-form $t_i^{-1}\rd t_i$.
Set $\mathcal K=\Omega^1_\mathcal R/\rd(\Omega_{\mathcal R}^1)$, where
\[\rd(\Omega_{\mathcal R}^1)=\text{Span}_{\C}\{r_0t_0^{r_0}t_1^{r_1}\rk_0+r_1t_0^{r_0}t_1^{r_1}\rk_1\mid r_0,r_1\in \Z\}\]
is the space of exact $1$-forms in $\Omega_\mathcal R^1$.
Let $\dg$ be a finite dimensional complex simple Lie algebra, and let $\la\cdot,\cdot \ra$ be the normalized non-degenerate
invariant symmetric bilinear form on it. The universal central extension of the loop algebra $\mathcal R\ot \dg$, called $2$-toroidal Lie algebra, can be realized as
$\(\mathcal R\ot \dg\)\oplus \mathcal K,$
and its Lie bracket is given by (see \cite{MRY},\cite{RM})
\[ [t_0^{m_0}t_1^{m_1}\ot x, t_0^{n_0}t_1^{n_1}\ot y]=t_0^{m_0+n_0}t_1^{m_1+n_1}\ot [x,y]+\la x,y\ra \sum_{a=0,1} m_at_0^{m_0+n_0}t_1^{m_1+n_1}\rk_a,\]
where $x,y\in \dg$, $m_0,n_0, m_1,n_1\in \Z$,  and $\mathcal K$ is central.

We denote by $\mathcal D=\te{Der}(\mathcal R)$ the Lie algebra of derivations over $\mathcal R$.
As a left $\mathcal R$-module, $\te{Der}(\mathcal R)$ has a basis $\rd_0, \rd_1$, where $\rd_i=t_i\frac{\partial}{\partial t_i}$, $i=0,1$.
 The elements from the Lie algebra $\mathcal D$ act naturally on $\mathcal R\ot \dg$
\[[t_0^{m_0}t_1^{m_1}\rd_i, t_0^{n_0}t_1^{n_1}\ot x]=n_it_0^{m_0+n_0}t_1^{m_1+n_1}\ot x,\]
and this action can be uniquely extended to the center $\mathcal K$ of the toroidal Lie algebra $\(\mathcal R\ot \dg\)\oplus \mathcal K$
\[[t_0^{m_0}t_1^{m_1}\rd_i, t_0^{n_0}t_1^{n_1}\rk_j]=n_it_0^{m_0+n_0}t_1^{m_1+n_1}\rk_j+\delta_{i,j}
\sum_{a=0,1} m_at_0^{m_0+n_0}t_1^{m_1+n_1}\rk_a.\]

With the above Lie brackets, $\(\mathcal R\ot \dg\)\oplus \mathcal K \oplus \mathcal D$ forms a Lie algebra, which is often called a full toroidal Lie algebra of rank 2. This full toroidal Lie algebra contains a very important Lie subalgebra
\[\wt\fg=\(\mathcal R\ot \dg\)\oplus \mathcal K\oplus \mathcal S,\]
where $\mathcal S$ is a subspace of $\mathcal D$ consisting of
divergence zero derivations (also called skew derivations). That is
\[\mathcal S=\{d\in \mathcal D\mid \<d, \rd(\Omega_{\mathcal R}^1)\>=0\},\]
where $\<\cdot ,\cdot \>:\mathcal D\times \Omega_{\mathcal R}^1\rightarrow \C$ is a bilinear form determined by
\begin{align*}\label{pairing}
\<t_0^{m_0}t_1^{m_1}\rd_i, t_0^{n_0}t_1^{n_1}\rk_j\>=
\delta_{m_0+n_0,0}\delta_{m_1+n_1,0}\delta_{i,j}.
\end{align*}
 It is well known \cite{BGK} that this Lie algebra $\wt\fg$ admits a non-degenerate invariant bilinear form. We note that  the subspace $\mathcal S$
is spanned by the degree zero derivations $\rd_0, \rd_1$ and the following skew derivations
\[ \rd(m_0,m_1)=m_0 t_0^{m_0}t_1^{m_1}\rd_1-m_1 t_0^{m_0}t_1^{m_1}\rd_0,\]
for $ m_0,m_1\in \Z$.
 We know (see \cite{B}) that the Lie algebra $\mathcal S$ admits an abelian extension over $\mathcal K$ with the following $2$-cocycle $\tau_\mu:\mathcal S\times \mathcal S\rightarrow \mathcal K$ such taht
\begin{align*}
\tau_\mu(\rd(m_0,m_1),\rd(n_0,n_1))=\mu(m_0n_1-m_1n_0)^{2}\sum_{a=0,1} m_at_0^{m_0+n_0}t_1^{m_1+n_1}\rk_a,
\end{align*}
and $\tau_\mu(\rd_i,\mathcal S)=0$, where  $i=0,1$, and $\mu$ is a fixed complex number. By twisting this $2$-cocycle to the Lie algebra $\wt\fg$, one  obtains an EALA of nullity $2$
\[\widetilde{\fg}(\mu)=\(\mathcal R\ot \dg\)\oplus \mathcal K\oplus_{\tau_\mu} \mathcal S,\]
and the remaining bracket relations on $\wt\fg(\mu)$ are given by
\begin{eqnarray*}
&&[\rd_i,\rd(m_0,m_1)]=m_i\rd(m_0,m_1),\ [\rd_1,\rd_2]=0,\\
&&[\rd(m_0,m_1),\rd(n_0,n_1)]=(m_0n_1-m_1n_0)\rd(m_0+n_0,m_1+n_1)\\
&&\qquad \qquad \qquad +\mu(m_0n_1-m_1n_0)^{2}\sum_{a=0,1} m_at_0^{m_0+n_0}t_1^{m_1+n_1}\rk_a,
\end{eqnarray*}
for $x\in \dg$, $m_0,m_1,n_0,n_1\in \Z$ and $i,j=0,1$.
The resulting Lie algebra $\wt\fg(\mu)$ is often called nullity $2$ toroidal EALA, and this is the extended affine Lie algebra that will be studied in this paper.

Let $\dfh$ be the Cartan subalgebra of $\dg$.
 Then the Cartan subalgebra of $\widetilde{\fg}(\mu)$ is
given by
\[\widetilde{\fh}=\dfh\oplus \C\rk_0\oplus \C\rk_1 \oplus \C\rd_0 \oplus \C\rd_1,\]
and the non-degenerate invariant symmetric bilinear form $\<\cdot , \cdot\>$ on $\widetilde{\fg}(\mu)$ can be determined as follows
\begin{eqnarray*}
&&\la t_0^{m_0}t_1^{m_1}\ot x, t_0^{n_0}t_1^{n_1}\ot y\ra=\delta_{m_0+n_0,0}\delta_{m_1+n_1,0}
\la x, y\ra,\ \la \rd_i, \rk_j \ra=\delta_{i,j},\\
&&\la  \rd(m_0,m_1),  t_0^{n_0}t_1^{n_1}\rk_j\ra=
(m_1\delta_{0,j}-m_0\delta_{1,j})\delta_{m_0+n_0,0}\delta_{m_1+n_1,0},
\end{eqnarray*}
and the others are trivial.

Let $\dot{\Delta}$ be the root system of $\dg$ relative to $\dfh$.
For each root $\dal\in \dot{\Delta}$, we extend $\dal$ to be a linear functional in $\wt\fh^*$ so that
\[\dal(\rk_i)=0=\dal(\rd_i),\quad i=0,1,\]
and define linear functionals $\delta_i$, $i=0,1$ on $\wt\fh$  by setting
\begin{align*}&\delta_i(\dfh)=0,\quad \delta_i(\rk_j)=0,\quad \delta_i(\rd_j)=\delta_{i,j}.
\end{align*}
Then the set
\[\wt\Delta=\{\dal+m_0\delta_0+m_1\delta_1, m_0\delta_0+m_1\delta_1\mid \dal\in \dot\Delta,\ m_0,m_1\in \Z\}\]
is the root system of $\wt\fg(\mu)$ with respect to the Cartan subalgebra $\wt\fh$. And we have the following root space decomposition
\[\wt\fg(\mu)=\bigoplus_{\al\in \wt\Delta}\wt\fg(\mu)_\al,\]
where $\wt\fg(\mu)_\al=\{x\in \wt\fg(\mu)\mid [h,x]=\al(h)x,\ \forall h\in \wt\fh\}$.

Note that the affine Kac-Moody algebra $\fg$ defined in \eqref{affalg} is a Lie subalgebra of $\wt\fg(\mu)$, and the set
\begin{align}\label{affroot}
\Delta=\{\dal+m_0\delta_0, m_0\delta_0\mid \dal\in \dot\Delta,\ m_0\in \Z\}\subset \wt\Delta\end{align}
is the root system of $\fg$ relative to its Cartan subalgebra $\fh$.
We fix a simple root system $\Pi$ in $\Delta$, and
denote by $\Delta_+$ and $\Delta_-$ respectively the corresponding positive and negative root systems of $\fg$.
Then we have a decomposition
\[\wt\Delta=\wt\Delta_+\cup \wt\Delta_0\cup \wt\Delta_-\]
of $\wt\Delta$, where
$\wt\Delta_\pm=\{\al+m_1\delta_1\mid \al\in \Delta_\pm, m_1\in \Z\}$
 and
 $\wt\Delta_0=\{m_1\delta_1\mid m_1\in \Z\}$.
This gives a triangular decomposition
\[\wt\fg(\mu)=\wt\fg(\mu)_+\oplus \wt{\mathcal H}\oplus \wt\fg(\mu)_-\]
for the  nullity $2$ toroidal EALA $\wt\fg(\mu)$, where  $\wt\fg(\mu)_\pm=\oplus_{\al\in \wt\Delta_\pm}\wt\fg(\mu)_\al$, and
\begin{align}\label{ch}
\wt{\mathcal H}=\oplus_{\al\in \wt\Delta_0}\wt\fg(\mu)_\al=\(\C[t_1,t_1^{-1}]\ot \dfh\)\oplus \sum_{n\in \Z}\(\C t_1^n\rk_0\oplus \C t_1^n\rd_0\)\oplus \C\rk_1\oplus \C\rd_1.\end{align}

\section{Representations of the Virasoro-affine algebras}

Recall that the Virasoro-affine algebra $\bar\ff$  is related to the reductive Lie algebra $\dot\ff=\dot\fg\oplus \dot\fb$ (see \eqref{viraff}), and its
 Lie bracket relations are given by
\begin{align*}
&[L(m), L(n)]=(m-n)L(m+n)+\frac{m^3-m}{12}\delta_{m+n,0}\rk_v,\\
&[x(m), y(n)]=[x,y](m+n)+\delta_{m,-n}m\<x,y\>\rk_0,\\
&[L(m), x(n)]=-nx(m+n),\quad
[\rk_0,\bar\ff]=0=[\rk_v,\bar\ff],
\end{align*}
where $m,n\in \Z$, $x,y\in \dot\ff$, $x(m)=t_0^m\ot x$ and $L(m)=-t_0^m\rd_0=-t_0^{m+1}\frac{d}{d t_0}$.
Recall also that $\bar\ff_0$ is a Cartan subalgebra of $\bar\ff$ (see \eqref{barh}).
We view the affine root system $\Delta$  (see \eqref{affroot}) as a subset of $\bar\ff_0^*$ by letting
\[\al(\rk_v)=\al(\rd_1)=\al(\rk_1)=0,\quad \al\in \Delta.\]
Then we have the following root space decomposition of $\bar\ff$
\[\bar\ff=\bigoplus_{\al\in \Delta}\bar\ff_\al,\]
where $\bar\ff_\al=\{x\in \bar\ff\mid [h,x]=\al(h)x,\ \forall h\in \bar\ff_0\}$.
This induces  a triangular decomposition of $\bar\ff$ as follows
\begin{align}\label{tridec}
\bar\ff=\bar\ff_+\oplus \bar\ff_0\oplus \bar\ff_-,\end{align}
where $\bar\ff_\pm=\oplus_{\al\in \Delta_\pm}\bar\ff_\al$.
For any given linear functional $\eta$ on $\bar\ff_0$, we denote by $\C_\eta=\C v_\eta$ the $1$-dimensional $(\bar\ff_+\oplus \bar\ff_0)$-module such that
 \[\bar\ff_+.v_\eta=0,\quad h.v_\eta=\eta(h)v_\eta,\ h\in \bar\ff_0.\]
Then we have the induced $\bar\ff$-module
\[M_{\bar\ff}(\eta)=\mathcal U(\bar\ff)\otimes_{\mathcal U(\bar\ff_+\oplus \bar\fh)}\C_\eta.\]
Here and henceforth, the notation $\mathcal U(\cdot)$ stands for the universal enveloping algebra of a given Lie algebra.
Write $V_{\bar\ff}(\eta)$ for the irreducible quotient of $M_{\bar\ff}(\eta)$, which is the irreducible highest weight $\bar\ff$-module
with highest weight $\eta$.

Note that the affine Kac-Moody algebra $\fg$ is a subalgebra of $\bar\ff$, and $\bar\ff$ contains also a Virasoro subalgebra $\mathfrak V=\mathrm{Der}(\C[t_0,t_0^{-1}])\oplus \C\rk_v$ and
a generalized Heisenberg subalgebra
$\fb=(\C[t_0,t_0^{-1}]\ot \dot{\fb})\oplus \C\rk_0\oplus \C\rd_0$. We denote  by $\bar\fg$ the Virasoro-affine subalgebra
$\mathfrak V+\fg$.
If $\fv$ is any one of these Lie subalgebras of $\bar\ff$,
then the decomposition \eqref{tridec} induces a triangular  decomposition of $\fv$ as follows
\begin{align*}
\fv=\fv_+\oplus \fv_0\oplus \fv_-,\quad \fv_\pm=\fv\cap \bar\ff_\pm,\ \fv_0=\fv\cap \bar\ff_0.
\end{align*}
Using this decomposition, for any $\eta\in \fv_0^*$, one can define the irreducible highest weight $\fv$-module $V_{\fv}(\eta)$ with highest weight $\eta$ as we did for $\bar\ff$.

From now on, let $\lambda$ be a fixed linear functional on $\fh$ such that $c=\lambda(\rk_0)\ne 0$.
 View the Laurent polynomial ring
 $\C[t_1,t_1^{-1}]$ as a $(\fb_+\oplus \fb_0)$-module by letting
\[x.t_1^n=\<x,n\rk_1\>t_1^n,\quad \rk_0.t_1^n=c\ t_1^n, \quad \rd_0.t_1^n=0, \quad \fb_+.t_1^n=0,\]
where $x\in \dot\fb$ and $n\in \Z$.
Then we have the following induced $\fb$-module
\[L_{\fb}(c)=\mathcal U(\fb)\ot_{\mathcal U(\fb_+\oplus \fb_0)}\C[t_1,t_1^{-1}]=\mathrm{S}(\fb_-)\ot \C[t_1,t_1^{-1}],\]
where $\mathrm{S}(\fb_-)$ stands for the symmetric algebra over $\fb_-$.
We extend $\lambda$ to be a linear functional $\tilde\lambda$ on $\bar\fg_0$ by letting
\[\tilde{\lambda}(\rk_v)=24\mu c-2,\]
and  define a $\bar\ff$-module structure on the tensor product space
\[L_{\bar\ff}(\lambda)=V_{\bar\fg}(\tilde{\lambda})\ot  L_{\fb}(c)\] as follows
\begin{align*}
&x(z)=\sum_{m\in \Z} x(m)z^{-m-1}\mapsto x(z)\ot 1,\quad x\in \dg,\\
&y(z)=\sum_{m\in \Z}y(n)z^{-m-1}\mapsto 1\ot y(z),\quad y\in \dot{\fb}, \\
& L(z)=\sum_{m\in \Z} L(m)z^{-m-2}\mapsto L(z)\ot 1+1\ot \frac{1}{c}:\rd_1(z)\rk_1(z):,\\
&\rk_0\mapsto c\,\mathrm{Id}, \quad \rk_v\mapsto 24\mu c\,\mathrm{Id},\end{align*}
where and hereafter, the normal ordered product $:a(z)b(w):$ for any two fields $a(z)=\sum_{n\in \Z}a_{(n)}z^{-n-1}$ and $b(w)$ over any vector space  is defined as follows
\[:a(z)b(w):=a(z)_-b(w)+b(w)a(z)_+,\]
where $a(z)_-=\sum_{n<0}a_{(n)}z^{-n-1}$ and $a(z)_+=\sum_{n\ge 0}a_{(n)}z^{-n-1}$. For each $n\in \Z$, we set
\begin{align}Y(n,z)=E^-(n,z)E^+(n,z)\, t_1^n,\end{align}
where $t_1^n$ is regarded as a left multiplication operator on $\C[t_1,t_1^{-1}]$, and
\begin{align}\label{enz}
E^\pm(n,z)=\mathrm{exp}\(\frac{1}{c}\sum_{ j\ge 1}\frac{n\rk_1(\pm j)}{\mp j}z^{\mp j}\).
\end{align}
The following result is given in \cite{B}, and which will be used in our paper to prove Proposition 3.2.
\begin{thm}\label{thm:b2} Let $\lambda$ be a linear functional on $\fh$ such that $c=\lambda(\rk_0)\ne 0$. Then there is a $\wt\fg(\mu)$-module structure on  the $\bar\ff$-module
 $L_{\bar\ff}(\lambda)$ with
 actions given by
\begin{align*}
&\sum_{m\in \Z}t_0^mt_1^n\rk_0 z^{-m}\mapsto c\,Y(n,z),\quad
\sum_{m\in \Z}t_0^mt_1^n\rk_1 z^{-m-1}\mapsto \rk_1(z)Y(n,z),\\
&\sum_{m\in \Z}t_0^mt_1^n\ot x z^{-m-1}\mapsto x(z)Y(n,z),\quad
 \rd_1\mapsto \rd_1(0),\quad \rd_0\mapsto \rd_0,\\
 & \sum_{m\in \Z}\hat\rd(m,n)z^{-m-2}\mapsto n:L(z)Y(n,z):+
  n^2(\mu-\frac{1}{c})\(\frac{\partial}{\partial z}\rk_1(z)\)Y(n,z)\\
&\qquad\qquad\qquad\qquad\quad-\(z^{-1}+\frac{\partial}{\partial z}\):\rd_1(z)Y(n,z):,\\
 \end{align*}
where $x\in \dg$, $n\in \Z$ and $\hat\rd(m,n)=\rd(m,n)+n\mu(m+1)t_0^mt_1^n\rk_0$.
\end{thm}

We denote by
\[\wh\fg(\mu)=(\mathcal R\ot \dg)\oplus \mathcal K\oplus \sum_{m,n\in \Z}\C\rd(m,n)\oplus \C\rd_0,\]
and $\wh\fg(\mu)$ is clearly a Lie subalgebra of $\wt\fg(\mu)$ and $\wt\fg(\mu)=\wh\fg(\mu)\oplus \C\rd_1$.
It is routine to check that the $\bar\ff$-submodule
\[V_{\bar\fg}(\tilde{\lambda})\ot \mathrm{S}(\fb_-)\] of $L_{\bar\ff}(\lambda)$ is isomorphic to the irreducible highest weight $V_{\bar\ff}(\bar\lambda)$ with highest weight $\bar\lambda$, where $\bar\lambda\in \bar\ff_0^*$  is
defined in \eqref{barlambda}.
We are going to show that there is a $\wh\fg(\mu)$-module structure on
this irreducible highest weight $\bar\ff$-module. For this purpose, we introduce four more fields on the $\bar\ff$-module $L_{\bar\ff}(\lambda)=V_{\bar\fg}(\tilde{\lambda})\ot \mathrm{S}(\fb_-)\ot
\C[t_1,t_1^{-1}]$ as follows
\begin{align}\label{barl}\bar{L}(z)=L(z)-\frac{1}{c}:\rd_1(z)\rk_1(z):,\ \ \ \bar\rk_1(z)=\(z^{-1}+\frac{\partial}{\partial z}\)\rk_1(z),\\
\label{bard}\bar\rd_1(z)=\(z^{-1}+\frac{\partial}{\partial z}\)\rd_1(z),\ \ \  E(n,z)=E^-(n,z)E^+(n,z).\end{align}
Note that the field $\bar{L}(z)$ acts  on $V_{\bar\fg}(\tilde{\lambda})$, and $\bar\rk_1(z), \bar\rd_1(z), E(n,z)$  act on
$\mathrm{S}(\fb_-)$. We now define the following operators, for $m,n\in \Z$ and $ x\in \dg $,
\[\rk_0(m,n),\ \rk_1(m,n),\ \bar\rk_1(m,n),\ \bar\rd_1(m,n),\ \bar{L}(m,n),\ \bar{\rd}(m,n),\ x(m,n)\]
acting on $V_{\bar\ff}(\bar\lambda)=V_{\bar\fg}(\tilde{\lambda})\ot \mathrm{S}(\fb_-)$ via formal power series expansions
\begin{align*}
\sum_{m\in \Z}\rk_0(m,n)z^{-m}&=\lambda(\rk_0) E(n,z),\
\sum_{m\in \Z} \rk_1(m,n)z^{-m-1}= \rk_1(z) E(n,z),\\
\sum_{m\in \Z} \bar\rk_1(m,n)z^{-m-2}&= \bar\rk_1(z) E(n,z),\
\sum_{m\in \Z} \bar\rd_1(m,n)z^{-m-2}= :\bar\rd_1(z) E(n,z):,\\
\sum_{m\in \Z}\bar{L}(m,n)z^{-m-2}&=\bar{L}(z) E(n,z),\
\sum_{m\in \Z} x(m,n)z^{-m-1}=x(z) E(n,z),\\
\bar{\rd}(m,n)&=n\bar{L}(m,n) +n^2\mu\bar\rk_1(m,n)-\bar\rd_1(m,n).
\end{align*}
\begin{prpt}\label{refor}
Assume that $\lambda\in \fh^*$ such that $c=\lambda(\rk_0)\ne 0$. Then there is a $\wh\fg(\mu)$-module structure on the
highest weight $\bar\ff$-module $V_{\bar\ff}(\bar\lambda)$ with actions given by
\begin{align*}
t_0^{m}t_1^{n}\rk_j\mapsto \rk_j(m,n),\ t_0^{m}t_1^{n}\ot x\mapsto x(m,n),\
\tilde\rd(m,n)\mapsto \bar{\rd}(m,n),\ \rd_0\mapsto \rd_0,
\end{align*}
where $m,n\in \Z, j=1,2$, $x\in \dg$ and $\tilde{\rd}(m,n):=\rd(m,n)+n\mu t_0^mt_1^n\rk_0$.
\end{prpt}
\begin{proof} Since $Y(n,z)=E(n,z)t_1^n$, and $\bar\rd_1(z)$ commutes with $t_1^n$, one has  $$:\bar\rd_1(z)Y(n,z):=:\bar\rd_1(z)E(n,z):t_1^n.$$ In view of this, by comparing the actions given respectively in Theorem \ref{thm:b2}, it suffices to show that
 the action of $\hat{\rd}(m,n)$ given in Theorem \ref{thm:b2} can be rewritten as follows
\begin{align}\label{rewrd}
\sum_{m\in \Z}\tilde{\rd}(m,n)z^{-m-2}\mapsto n\bar{L}(z)Y(n,z)
+\mu n^2 \bar{\rk}_1(z)Y(n,z)-:\bar{\rd}_1(z)Y(n,z):.
\end{align}
Note that $\frac{\partial}{\partial z}$ is a derivation of the normal ordered product (\cite[(3.1.5)]{K})
and $\frac{\partial}{\partial z}Y(n,z)=\frac{n}{c}\rk_1(z)Y(n,z)$. Therefore we have
\begin{equation}\begin{split}\label{rew1}
&\(z^{-1}+\frac{\partial}{\partial z}\):\rd_1(z)Y(n,z):\\
=\,&:\(\(z^{-1}+\frac{\partial}{\partial z}\)\rd_1(z)\)Y(n,z):
+:\rd_1(z) \(\frac{\partial}{\partial z}Y(n,z)\):\\
=\,&:\(\(z^{-1}+\frac{\partial}{\partial z}\)\rd_1(z)\)Y(n,z):
+\frac{n}{c}:\rd_1(z) :\rk_1(z)Y(n,z)::.
\end{split}\end{equation}
The following OPEs are well-known  (\cite[(5.4.3a)]{K})
\begin{align*}
[\rd_1(z),Y(n,w)]=nw^{-1}\delta(w/z)Y(n,w),\quad [\rk_1(z),Y(n,w)]=0,
\end{align*}
where the $\delta$-function $\delta(w/z)=\sum_{m\in \Z}(w/z)^m$.
Using these OPEs and the ``quasiassociativity" of the normal ordered product (\cite[(4.8.5)]{K}), it is easy to check that
\begin{align}\label{rew2}
::\rd_1(z)\rk_1(z):Y(n,z):-:\rd_1(z) :\rk_1(z)Y(n,z)::=
n\(\frac{\partial}{\partial z}\rk_1(z)\)Y(n,z).
\end{align}
Moreover,  as fields acting on $L_{\bar\ff}(\lambda)$, one has
\begin{align}\label{rew3}
\sum_{m\in \Z}\tilde{\rd}(m,n)z^{-m-2}-\sum_{m\in \Z}\hat{\rd}(m,n)z^{-m-2}=\mu n^2(z^{-1}\rk_1(z))Y(n,z).
\end{align}
It is then easy to see that \eqref{rewrd} follows from \eqref{barl}, \eqref{bard}, \eqref{rew1}, \eqref{rew2} and \eqref{rew3}.
This completes the proof of the proposition.
\end{proof}

To end this section, we recall the following lemma for later use, and the proof this lemma follows from the well-known Segal-Sugawara construction \cite[Proposition 3.1]{B}.
\begin{lemt}\label{hwmdec}
Let $\mathrm{h}^\vee$ be the dual Coxeter number of $\dg$, and $c\ne -\mathrm{h}^\vee$.  Then
one has the following decomposition of the $\bar\fg$-module $V_{\bar\fg}(\tilde{\lambda})$
\[V_{\bar\fg}(\tilde{\lambda})=V_{\fg}(\lambda)\ot V_{\mathfrak V}(\lambda_{\mathfrak V})\]
where $\lambda_{\mathfrak V}\in \mathfrak{V}_0^*$ is defined  by
\[\lambda_{\mathfrak V}(\rk_v)=24\mu c-2-\frac{\eta(\rk_0)\dim(\dg)}{c+\mathrm{h}^\vee},\quad
 \lambda_{\mathfrak V}(\rd_0)=\lambda(\rd_0)-\frac{\<\lambda|_{\dfh},\lambda|_{\dfh}\>}{2(c+\mathrm{h}^\vee)}.\]
\end{lemt}

\section{Loop representations of $\wt\fg(\mu)$}
In this section we construct a class of modules $\wt V(\bm{\lambda},\bm{a})$ for the nullity 2 toroidal extended affine Lie algebra $\wt\fg(\mu)$. The reducibility and integrability of the $\wt\fg(\mu)$-module $\wt V(\bm{\lambda},\bm{a})$ will be given respectively in Theorem 4.5 and Proposition 4.7.

Note that the algebra $\wt\fg(\mu)$ is $\Z$-graded with respect to the action of $\rd_1$, and
\[\wt\fg(\mu)=\oplus_{n\in \Z}\wt\fg(\mu)_n,\]
where $\wt\fg(\mu)_n=\{x\in \wt\fg(\mu)\mid [\rd_1,x]=nx\}$. With this gradation, if $\mathcal G$ is any graded subalgebra of $\wt\fg(\mu)$,
 then there is a natural $\Z$-grading on its universal enveloping algebra $\U(\mathcal G)$, and the corresponding homogenous subspaces will be denoted respectively by $\mathcal G_n$ and $\U(\mathcal G)_n$ for $n\in \Z$.

Let $k$ be a fixed positive integer, and $(\bm{\lambda},\bm{a})\in (\fh^*)^k\times (\C^\times)^k$ be a fixed pair as in \eqref{lambdaa} and satisfy the condition \eqref{assumption}.
By  Proposition \ref{refor}, we know that
there is a $\wh\fg(\mu)$-module structure on the tensor product space
\[\wh V(\bm{\lambda},\bm{a})=V_{\bar\ff}(\bar\lambda_1)\ot\cdots \ot V_{\bar\ff}(\bar\lambda_k)\]
with actions given by
\begin{align}
\label{act1}&t_0^{m}t_1^{n}\rk_j. (v_1\ot \cdots \ot v_k)=\sum_{i=1}^k a_i^n v_1\ot \cdots \ot \rk_j(m,n).v_i\ot \cdots \ot v_k,\\
\label{act2}&t_0^{m}t_1^{n}\ot x. (v_1\ot \cdots \ot v_k)=\sum_{i=1}^k a_i^n v_1\ot \cdots \ot x(m,n).v_i\ot \cdots \ot v_k,\\
\label{act3}&\tilde\rd(m,n). (v_1\ot \cdots \ot v_k)=\sum_{i=1}^k a_i^n v_1\ot \cdots \ot \bar\rd(m,n).v_i\ot \cdots \ot v_k,\\
&\rd_0. (v_1\ot \cdots \ot v_k)=\sum_{i=1}^k  v_1\ot \cdots \ot \rd_0.v_i\ot \cdots \ot v_k,
\end{align}
where $m,n\in \Z, j=0,1,$ and $ x\in \dg$, $v_i\in V_{\bar\ff}(\bar\lambda_i)$, $i=1,\cdots,k$. Similar to the loop module construction for affine Kac-Moody algebras, there is a $\wt\fg(\mu)$-module structure on the loop space
\[\wt V(\bm{\lambda},\bm{a})=\wh V(\bm{\lambda},\bm{a})\ot \C[t_1,t_1^{-1}]\]
with the actions given by
\begin{align}
x.(v\ot t_1^l)=(x.v)\ot t_1^{l+n},\quad \rd_1.(v\ot t_1^l)=l\,v\ot t_1^l,
\end{align}
where $x\in \wh\fg(\mu)_n$, $v\in \wh V(\bm{\lambda},\bm{a})$ and $l\in \Z$.

We remark that the most crucial step in the study of the structure of the $\wt\fg(\mu)$-module $\wt V(\bm{\lambda},\bm{a})$ is the verification of the irreducibility of the $\wh\fg(\mu)$-module $\wh V(\bm{\lambda},\bm{a})$. It is easy to see that $\wh V(\bm{\lambda},\bm{a})$ is a weight module with respect to $\wh\fh$, and
\begin{align*}
\wh V(\bm{\lambda},\bm{a})=\oplus_{\gamma\in \wh\fh^*} \wh V(\bm{\lambda},\bm{a})_{\gamma},\end{align*}
where $\wh\fh=\wt\fh\cap \wh\fg$, and
$\wh V(\bm{\lambda},\bm{a})_{\gamma}=\{v\in V(\bm{\lambda},\bm{a})\mid h.v=\gamma(h)v,\forall\,h\in \wh\fh\}$.
Denote by $v_{\bar\lambda_i}$ a fixed highest weight vector of $V_{\bar\ff}(\bar\lambda_i)$, and set
\[v_{\bm{\lambda}}=v_{\bar\lambda_1}\ot \cdots \ot v_{\bar\lambda_k},\]
we define a linear functional $\wh\psi_{\bm{\lambda},\bm{a}}$ on  $\wh\CH=\wt\CH\cap \wh\fg(\mu)$ (see \eqref{ch}) by letting
\begin{align*}
t_1^n\ot \dot{h}\mapsto \sum_{i=1}^k a_i^n \lambda_i(\dot{h}),\quad
t_1^n \rk_0\mapsto \sum_{i=1}^k a_i^n \lambda_i(\rk_0),\quad
\rk_1\mapsto 0,\\
 \rd_0\mapsto \sum_{i=1}^k\lambda_i(\rd_0),\quad  t_1^m \rd_0\mapsto\sum_{i=1}^k a_i^m (\lambda_i(\rd_0)+\mu\lambda_i(\rk_0)),
\end{align*}
where  $\dot{h}\in \dfh$, $n\in \Z$ and $m\in {\Z\setminus \{0\}}$.
Then we have the following lemma. Its proof is straightforward and the argument is omitted for shortness.
\begin{lemt}\label{basic} (i). The weights of $\wh V(\bm{\lambda},\bm{a})$ have the form $\underline{\lambda}-\eta$ for some
 $\eta\in \Gamma_+$, where $\underline{\lambda}=\sum_{i=1}^k \lambda_i$ and $\Gamma_+=\oplus_{\alpha\in\Pi}{\N}\alpha$.

(ii). $\wt\fg(\mu)_+.v_{\bm{\lambda}}=0$ and $\wh V(\bm{\lambda},\bm{a})_{\underline{\lambda}}=\C v_{\bm{\lambda}}$.

(iii). $h.v_{\bm{\lambda}}=\wh\psi_{\bm{\lambda},\bm{a}}(h) v_{\bm{\lambda}}$ for all $h\in \wh\CH$.
\end{lemt}

By applying Lemma 4.1, one may prove the following key theorem. Due to the complexity and lengthy in the proof of the theorem, we would like to provide its proof in Section 5.
\begin{thm}\label{irr}
The $\wh\fg(\mu)$-module $\wh V(\bm{\lambda},\bm{a})$ is irreducible.
\end{thm}

Suppose that the $\wh\fg(\mu)$-module $\wh V(\bm{\lambda},\bm{a})$ is irreducible, we prove that the $\wt\fg(\mu)$-module $\wt V(\bm{\lambda},\bm{a})$ is completely reducible with finitely many irreducible components.
We first note that $\wt V(\bm{\lambda},\bm{a})$ is a weight module for $\wt\fg(\mu)$ with respect to the Cartan subalgebra $\wt\fh$, and all weights have the form
$\underline{\lambda}-\eta+l\delta_1$ for some $\eta\in \Gamma_+$, $l\in \Z$.
Set
$
 \Omega_{\bm{\lambda},l}=v_{\bm{\lambda}}\ot t_1^l,\ l\in \Z.
$
It is clear that $\wt V(\bm{\lambda},\bm{a})_{\underline{\lambda}+l\delta_1}=\C \Omega_{\bm{\lambda},l}$.
Define $\Z$-graded algebra homomorphism
\[\wt\psi_{\bm{\lambda},\bm{a}}:\mathcal U(\wh\CH)\rightarrow \C[t_1,t_1^{-1}],\]
by setting
 $
\wt\psi_{\bm{\lambda},\bm{a}}(h)=\wh\psi_{\bm{\lambda},\bm{a}}(h) t_1^n$ for $ h\in \wh\CH_n.$
Then, by Lemma \ref{basic} (iii), we have
\begin{align}\label{wtpsi}
h.\Omega_{\bm{\lambda},l}=\wt\psi_{\bm{\lambda},\bm{a}}(h)\Omega_{\bm{\lambda},n+l},\quad h\in \U(\wh\CH)_n.
\end{align}

The image of $\wt\psi_{\bm{\lambda},\bm{a}}$ is equal to $L_r=\C[t_1^r,t_1^{-r}]$ for some integer $r\ge 1$ (see \cite{C}). Moreover, by a similar argument as in \cite[Lemma 4.4]{CP}, one has the following result.
\begin{lemt}\label{charchi}  Assume that the image of $\wt\psi_{\bm{\lambda},\bm{a}}$ is equal to $L_r=\C[t_1^r,t_1^{-r}]$ for some integer $r\ge 1$. Then we have $k\equiv 0\,(\,\mathrm{mod}\ r)$. Moreover,
there exist a permutation $\tau$ of $\{1,2,\cdots,k\}$, and complex numbers $a_{(0)},\cdots, a_{(p-1)}$
 such that
\[\lambda_{\tau(sr+1)}=\lambda_{\tau(sr+2)}=\cdots =\lambda_{\tau((s+1)r)},\]
 and
\[a_{\tau(sr+1)}=\varepsilon a_{(s)},\ a_{\tau(sr+2)}=\varepsilon^2 a_{(s)}, \cdots, a_{\tau((s+1)r)}=\varepsilon^r a_{(s)},\]
for $s=0,1,\cdots,p-1,$
where $p=k/r$, and $\varepsilon$ is a primitive $r$-th root of unity.
\end{lemt}

As a consequence of Lemma \ref{charchi}, we may assume that
\[\lambda_1=\cdots=\lambda_r,\ \lambda_{r+1}=\cdots=\lambda_{2r},\cdots,
\lambda_{k-r+1}=\cdots=\lambda_k.\]
Let $\sigma$ be the following permutation
\[\sigma=(1,2,\cdots,r)(r+1,r+2,\cdots,2r)\cdots(k-r+1,k-r+2,\cdots,k).\]
Then it induces an automorphism of $\wt V(\bm{\lambda},\bm{a})$ as a vector space, still denoted as $\sigma$, with the action given by
\[\sigma(v_1\ot \cdots \ot v_k\ot t_1^n)=\varepsilon^n(v_{\sigma(1)}\ot \cdots \ot v_{\sigma(k)}\ot t_1^n),\]
for $
v_i\in V_{\bar\ff}(\bar\lambda_i),\ n\in \Z.$
One can easily check that $ x.\sigma(v)=\sigma(x.v)$ for $ x\in \wt\fg(\mu),\ v\in \wt V(\bm{\lambda},\bm{a}).$
Thus one has the
following decomposation
\be\label{comred}\wt V(\bm{\lambda},\bm{a})=\bigoplus_{i=0}^{r-1}\wt V^i(\bm{\lambda},\bm{a})\ee
as $\wt\fg(\mu)$-submodules,
where
\begin{align*}
\wt V^i(\bm{\lambda},\bm{a})=\{v\in V(\bm{\lambda},\bm{a})\mid \sigma.v=\varepsilon^{-i}v\},\end{align*}
for $ i=0,1,\cdots,r-1.$
The following result is obvious.
\begin{lemt}\label{charhwv}
For each $l\in \Z$ and $i=0,\cdots,r-1$, the weight vector $\Omega_{\bm{\lambda},l}\in \wt V^i(\bm{\lambda},\bm{a})$ if and only if
$l\equiv i\,(\,\te{mod}\ r)$.
\end{lemt}

 Now we prove the first main result of the paper.

\begin{thm}\label{com} The $\wt\fg(\mu)$-module $\wt V(\bm{\lambda},\bm{a})$ is completely reducible, and each component $\wt V^i(\bm{\lambda},\bm{a})$ in \hbox{\eqref{comred}} is an irreducible $\wt\fg(\mu)$-submodule for $i=0,1,\cdots, r-1$.
\end{thm}
\begin{proof} Denote by $W$ the $\wt\fg(\mu)$-submodule of $\wt V(\bm{\lambda},\bm{a})$ generated by the weight vectors $\Omega_{\bm{\lambda},i}$, $i=0,1,\cdots, r-1$. We want to show that $W=\wt V(\bm{\lambda},\bm{a})$. Let $v\ot t_1^l$ be a given vector in $\wt V(\bm{\lambda},\bm{a})$,
where $v\in \wh V(\bm{\lambda},\bm{a})$ and $l\in \Z$.
By using Theorem \ref{irr}, there exists an element $x\in \U(\wh\fg(\mu))$ such that $x.v_{\bm{\lambda}}=v$.
Write $x=\sum_{j=1}^t x_j$ with $x_j\in \U(\wh\fg(\mu))_{n_j}$ for some $n_j\in \Z$.
Then we have
\be\label{cyclic}\sum_{j=1}^t x_j. \Omega_{\bm{\lambda},l-n_j}=\sum_{j=1}^t (x_j.v_{\bm{\lambda}})\ot t_1^l=
v\ot t_1^l.\ee
Since the image of $\wt\psi_{\bm{\lambda},\bm{a}}$ is equal to $L_r$, it follows from \eqref{wtpsi} that
\be \label{wtch}\Omega_{\bm{\lambda},i+rl}\in \U(\wh{\mathcal H})\Omega_{\bm{\lambda},i},\ee
for all $ i,l\in \Z$. This, together with \eqref{cyclic}, gives  $\wt V(\bm{\lambda},\bm{a})=W$.
Therefore, one concludes from  \eqref{comred} and Lemma \ref{charhwv} that each $\wt\fg(\mu)$-submodule $\wt V^i(\bm{\lambda},\bm{a})$
is generated by $\Omega_{\bm{\lambda},i}$.

Furthermore, we suppose that $U$ is a non-zero $\wt\fg(\mu)$-submodule of $\wt V^i(\bm{\lambda},\bm{a})$.
Define a partial order ``$\preceq$" on the set $\{\underline{\lambda}-\eta+l\delta_1\mid \eta\in \Gamma_+, l\in \Z\}$
by letting
$\underline{\lambda}-\eta+l\delta_1 \preceq \underline{\lambda}-\eta'+l'\delta_1$ if and only if $ \eta-\eta'\in \Gamma_+.$
  Take a non-zero weight vector $u\ot t_1^l$  in $U$ so that its weight is maximal with respect to this partial order.
 Then it is easy to see that $\wt\fg(\mu)_+.u\ot t_1^l=0$.
 This implies that $\wt\fg(\mu)_+.u=0$, and hence $u\in \C v_{\bm\lambda}$ (see Proposition \ref{dhwv}).
Thus, by Lemma \ref{charhwv}, one gets that $\Omega_{\bm{\lambda},i+rj}\in U$ for some $j\in \Z$.
This, together with \eqref{wtch}, gives that $\Omega_{\bm{\lambda},i}\in U$ and hence
$U=\wt V^i(\bm{\lambda},\bm{a})$ as required.
\end{proof}

Next, we consider the integrability of the $\wt\fg(\mu)$-module $\wt V(\bm{\lambda},\bm{a})$.
By definition, a $\wt\fg(\mu)$-module, or a $\wh\fg(\mu)$-module, is called integrable if it is a weight module and for every $\al\in \wt\Delta^\times$,
$\wt\fg(\mu)_\al$ acts locally nilpotent on the module, where
\[\wt\Delta^\times=\{\dot\al+m_0\delta_0+m_1\delta_1\mid \al\in \dot\Delta, m_0,m_1\in \Z\}\subset \wt\Delta.\]
For each $\al\in \Delta$, we denote by $\al^\vee\in \fh$ the coroot of $\al$, and by
\[P_+=\{\lambda\in \fh^*\mid \lambda(\al^\vee)\in \N,\ \forall\, \al\in \Pi\}\]
 the set of dominant weights of the affine Kac-Moody algebra $\fg$.

\begin{lemt}\label{integ} Assume that $\lambda\in P_+\setminus\{0\}$. Then the $\wh\fg(\mu)$-module $V_{\bar\ff}(\bar\lambda)$ defined in Proposition \ref{refor} is integrable.
\end{lemt}
\begin{proof} Recall from Lemma \ref{hwmdec} we have $V_{\bar\ff}(\bar\lambda)=V_\fg(\lambda)\ot V_{\mathfrak V}(\lambda_{\mathfrak V})\otimes \mathrm{S}(\fb_-)$. Let $\al\in \wt\Delta^\times$, $t_0^mt_1^n\ot x\in \wt\fg_\al$,  and $v\ot w\ot u\in V_\fg(\lambda)\ot V_{\mathfrak V}(\lambda_{\mathfrak V})\ot \mathrm{S}(\fb_-)$.
 Since the highest weight $\fg$-module $V_\fg(\lambda)$ is integrable, there exists a
sufficiently large integer $N$ such that
\begin{align*}
x(m_1)\cdots x(m_N).v=0,\end{align*}
for all $ m_1,\cdots,m_N\in \Z$, and therefore
\[(t_0^mt_1^n\ot x)^N.(v\ot w\ot u)=x(m,n)^N.(v\ot w\ot u)=0.\]
This proves the integrability of the $\wh\fg(\mu)$-module $V_{\bar\ff}(\bar\lambda)$.
\end{proof}

\begin{prpt} Assume that $\bm{\lambda}\in (P_+)^k$. Then the  $\wt\fg(\mu)$-module
 $\wt V(\bm{\lambda},\bm{a})$ is integrable.
\end{prpt}
\begin{proof}
It is suffice to show that the $\wh\fg(\mu)$-module $\wh V(\bm{\lambda},\bm{a})$ is integrable.
 For $v_1\ot \cdots \ot v_k\in \wh V(\bm{\lambda},\bm{a})$ and $x\in \wt\fg(\mu)_\al, \al\in \wt\Delta^\times$,
by Lemma \ref{integ}, there exist positive integers $N_i$, for $i=1,\cdots,k$, such that
$x^{N_i}.v_i=0$. Set $N=\sum_{i=1}^k N_i$. Then it is clear that
$x^N.(v_1\ot \cdots \ot v_k)=0$. This proves the integrability of the $\wh\fg(\mu)$-module $\wh V(\bm{\lambda},\bm{a})$.  Therefore, the $\wt\fg(\mu)$-module
$\wt V(\bm{\lambda},\bm{a})$ is also integrable.
\end{proof}

\section{Proof of Theorem \ref{irr}}

This section is devoted to a proof of Theorem \ref{irr}, and which follows from a sequence of lemmas.
For the sake of convenience, in this section we set $I=\{1,2,\cdots,k\}$ and $c_i=\lambda_i(\rk_0)\in \C^\times$, $i\in I$.
 The following lemma is given in \cite[Lemma 2.1]{BZ}, which is indeed a consequence of Vandermonde type determinant.

\begin{lemt}\label{van}Let $N$ be a non-negative integer. Then the following $k(N+1)\times k(N+1)$-matrix
\[(a_{ij})_{1\le i,j\le k(N+1)}\] is invertible, where
\[a_{ij}=a_{\underline{i}}^j\cdot j^{\bar{i}},\quad 1\le i,j\le k(N+1)\] and
$\bar{i}, \underline{i}$ are non-negative integers determined by
\[i=\bar{i}k+\underline{i},\quad 0\le \bar{i}\le N,\quad 1\le \underline{i}\le k.\]
\end{lemt}

Recall that the operators $E^\pm(n,z)$,  defined in \eqref{enz}, acting on  $V_{\bar\ff}(\bar\lambda_i)$ for $i\in I$.  We set
\begin{align*}
E^\pm(n,z)=\sum_{m\ge 0}\phi_n^\pm(m)z^{\mp m},
\end{align*}
and
\begin{equation}\phi^\pm(m,s)=\begin{cases}\delta_{m,0},\ &\te{if}\ s=0;\\
\frac{1}{s!}\sum_{m_1,\cdots,m_s>0; m_1+\cdots+m_s=m}\rk_1^\pm(m_1)\cdots \rk_1^\pm(m_s)\ &\text{if}\ s>0,\end{cases}\end{equation}
for $m\in \N$ and $0\le s\le m$, where $\rk_1^\pm(m_i):=\mp\frac{1}{\lambda_i(\rk_0)}\frac{\rk_1(\pm m_i)}{m_i}$.
Then it is easy to see that
\begin{align}\label{deck1}
\phi_n^\pm(m)=\sum_{0\le s\le m} n^s \phi^\pm(m,s),
\end{align}
for $  m\in \N,\ n\in \Z$.

Suppose that $W$ is a non-zero $\wh\fg(\mu)$-submodule of $\wh V(\bm{\lambda},\bm{a})$. Let
 $\bm{v}$ be a given non-zero element of $W$ with the form
\[\bm{v}=v_1\ot v_2\ot \cdots \ot v_k,\]
where $v_i\in V_{\bar\ff}(\bar\lambda_i)$. Note that
 there is a natural $\Z$-grading on $V_{\bar\ff}(\bar\lambda_i)$ with respect to the action of $\rd_0$. Namely,
\begin{align*}
V_{\bar\ff}(\bar\lambda_i)=\oplus_{n\in \N }V_{\bar\ff}(\bar\lambda_i)_{n},\end{align*}
where $V_{\bar\ff}(\bar\lambda_i)_{n}=
\{w\in V_{\bar\ff}(\bar\lambda_i)\mid \rd_0.w=(\lambda_i(\rd_0)-n)w\}.$
Fix a sufficiently large integer $M$ so that
\be\label{vanish}v_i\in \bigoplus_{n\le M}V_{\bar\ff}(\bar\lambda_i)_n\ee
for $ i\in I$.

\begin{lemt}\label{lem1} For $m\in \Z$ and $i\in I$, one has
\begin{align*}
v_1\ot \cdots \ot\rk_1(m).v_i\ot \cdots \ot v_k\in W.
\end{align*}
\end{lemt}
\begin{proof} From \eqref{vanish} we see that $\rk_1(m).v_i=0$ for all $i\in I$ and $m>M$.
Therefore, we may assume that $m\le M$.
For each $0\le s\le 2M-m$, we set
\be \label{L21}\phi^M(m,s)=\sum \phi^-(m_1,s_1)\phi^+(m_2,s_2),\ee
where the summation  is taken over the set
\[\{(m_1,m_2,s_1,s_2)\in \N^4\mid  m_2\le M, m_1=m_2-m, s_1\le m_1, s_2\le m_2, s_1+s_2=s\},\]
and also set
\begin{equation}\label{L22}
\phi_n^M(m)=
\sum_{0\le s\le 2M-m}n^s \phi^M(m,s),
\end{equation}
for  $n \in \Z$.
It follows from \eqref{vanish} that
\[\phi_n^+(m).v_i=0,\]
for $ n\in \Z, i\in I, m>M.$
Using this and \eqref{deck1}, it is easy to see that
\be\label{L23}\rk_0(m,n).v_i=c_i\,\phi^M_n(m).v_i,\ee
and
\be\label{L24} \rk_1(m).v_i=-mc_i\phi^M(m,1).v_i,\ee
for $ i\in I.$
Set $N=2M-m$. By  \eqref{act1}, we have the following  system of $k(N+1)$ equations
\begin{align*}
 t_0^{m}t_1^n\rk_0.\bm{v}=\sum_{i=1}^k a_i^n v_1\ot\cdots \ot \rk_0(m,n).v_i\ot \cdots\ot v_k,
\end{align*}
for $n=1,\cdots,k(N+1).$
 The above system of equations can be rewritten as follows by applying \eqref{L23},
\bee t_0^{m}t_1^n\rk_0.\bm{v}=\sum_{i=1}^k c_i a_i^n v_1\ot \cdots \ot \phi_n^M(m).v_i\ot \cdots\ot v_k,
\eee
for $ n=1,\cdots,k(N+1).$
This, together with \eqref{L22}, gives
\be\label{quo1} t_0^{m}t_1^n\rk_0.\bm{v}=\sum_{i=1}^k\sum_{s=0}^N  a_{i}^n n^s (c_i v_1\ot \cdots \ot \phi^{M}(m,s).v_i\ot \cdots \ot v_k),
\ee
for $n=1,\cdots,k(N+1)$.
Due to Lemma \ref{van}, the coefficient matrix of the above system of equations is invertible.
Thus, for every $i\in I$ and  $0\le s\le N$,
the vector
\be \label{phims}  c_i\, v_1\ot \cdots \ot \phi^{M}(m,s).v_i\ot \cdots \ot v_k\ee
is a linear combination of the elements
 $t_0^{m}t_1^n\rk_0.\bm{v}$ for $ n=1,2,\cdots,k(N+1).$ Therefore, the result of the lemma follows from \eqref{L24} and  \eqref{phims} by taking $s=1$.
\end{proof}

\begin{lemt}\label{lem2} For $x\in \dot{\fg}$, $m\in \Z$ and $i\in I$, one has
\begin{align*}
v_1\ot \cdots \ot x(m).v_i\ot \cdots \ot v_k\in W.
\end{align*}

\end{lemt}
\begin{proof} The argument is similar to that given in Lemma \ref{lem1}.
 By \eqref{vanish}, we may assume that $m\le M$.
For each $0\le s\le 3M-m$ and $n\in \Z$, we set
\begin{align*}
x^M(m,s)&=\sum_{m_1,m_2\le M, m_1+m_2=m} x(m_1)\phi^{M}(m_2,s),
\end{align*}
where $\phi^M(m_2,s)$ is defined in \eqref{L21}.
Then it is easy to check that
\[x(m,n).v_i=\sum_{0\le s\le 3M-m}n^s x^M(m,s).v_i,\]
for $ i\in I.$
This and \eqref{act2} give us the following system of equations for $1\le n\le k(3M-m+1)$
\begin{align}\label{quo2}
(t_0^mt_1^n\ot x).\bm{v}
=\sum_{i=1}^k \sum_{s=0}^{3M-m} a_i^n n^s (v_1\ot \cdots \ot x^M(m,s).v_i\ot \cdots \ot v_k).\end{align}
By a similar argument as we did in the proof of Lemma \ref{lem1}, one can solve this system of linear equations by applying Lemma \ref{van} to obtain
\[v_1\ot \cdots \ot x^M(m,s).v_i\ot \cdots \ot v_k\in W,\]
for $0\le s\le 3M-m$ and $i\in I$.
Then the lemma follows from this and the fact that $x^M(m,0).v_i=x(m).v_i$ for all $i\in I$.
\end{proof}

\begin{lemt}\label{lem3} For $m\in \Z$ and $i\in I$, one has
\begin{align*}
v_1\ot \cdots \ot \rd_1(m).v_i\ot \cdots \ot v_k,\
v_1\ot \cdots \ot L(m).v_i\ot \cdots \ot v_k\in W.
\end{align*}
\end{lemt}
\begin{proof}  For each $0\le s\le 3M-m$ and $n\in \Z$, we set
\begin{align*}
\bar{L}^M(m,s)&=\sum_{m_1,m_2\le M, m_1+m_2=m} \bar{L}(m_1)\phi^{M}(m_2,s),\\
\bar{\rk}_1^M(m,s)&=\sum_{m_1,m_2\le M, m_1+m_2=m} -m_1\rk_1(m_1)\phi^{M}(m_2,s),\\
\bar{\rd}_1^M(m,s)&=\sum_{m_1,m_2\le M, m_1+m_2=m} -m_1:\rd_1(m_1)\phi^{M}(m_2,s):,
\end{align*}
where $\phi^M(m_2,s)$ is defined in \eqref{L21}, and
\begin{equation}\label{norpro}
:\rd_1(m_1)\phi^{M}(m_2,s):=\begin{cases}
\rd_1(m_1)\phi^{M}(m_2,s), \quad &\te{if}\quad m_1\le 0;\\
\phi^{M}(m_2,s)\rd_1(m_1), \quad &\te{if}\quad m_1> 0.
\end{cases}\end{equation}
By applying \eqref{vanish}, one obtains
\begin{align*}
\rk_1(m).v_i=\rd_1(m).v_i=L(m).v_i=0,
\end{align*}
for $ i\in I, m>M, n\in \Z.$ Therefore,
we may assume that $m\le M$. Moreover, it follows from this and \eqref{L23} that
\begin{align*}
\bar{A}(m,n).v_i=\sum_{0\le s\le 3M-m}n^s \bar{A}^M(m,s).v_i,
\end{align*}
for $ i\in I,$
where $A=L$, $\rk_1$ or $\rd_1$.
This implies that
\begin{align}\label{bardmn}
\bar{\rd}(m,n).v_i=\sum_{0\le s\le 3M-m+2}n^s \bar{\rd}^M(m,s).v_i,\quad i\in I,
\end{align}
where $\bar{\rd}(m,n)$ is defined in Section 3, and
\begin{align*}
\bar{\rd}^M(m,s)=\bar{L}^M(m,s-1)+\mu\bar\rk_1^M(m,s-2)-\bar{\rd}_1^M(m,s),
\end{align*}
for $0\le s\le 3M-m+2.$
Here, as a convention we assume that $\bar{A}^M(m,s)=0$ for $s\notin \{0,1,\cdots,3M-m\}$. Therefore, from  \eqref{act3} and \eqref{bardmn}, we have the following system of equations for $1\le n\le 3M-m+3$
\begin{align}\label{quo3}
\tilde{\rd}(m,n).\bm{v}=\sum_{i=1}^k\sum_{s=0}^{3M-m+2} a_i^n n^s(v_1\ot \cdots \ot \bar{\rd}^M(m,s).v_i\ot \cdots \ot v_k),
\end{align}
and similarly as we did in the previous lemmas, we obtain  by applying Lemma \ref{van} that
\begin{align}\label{dms}
v_1\ot \cdots \ot \bar{\rd}^M(m,s).v_i\ot \cdots \ot v_k\in W,
\end{align}
for all $i\in I$ and $0\le s\le 3M-m+2$. Moreover, we note that $\rd_1(0)$ acts trivially on $V_{\bar\ff}(\bar\lambda_i)$, and
\begin{align*}
\bar{\rd}^M(m,0).v_i=-\bar{\rd}_1^M(m,0).v_i=m\rd_1(m).v_i,
\end{align*}
for $m\in \Z, i\in I.$ Then, by taking $s=0$ in \eqref{dms}, we get
\begin{align*}
v_1\ot \cdots \ot \rd_1(m).v_i\ot \cdots \ot v_k\in W,
\end{align*}
for $m\in \Z, i\in I.$ This, together with Lemma \ref{lem1}, gives that
\begin{align}\label{rkrd}
v_1\ot \cdots \ot :\rd_1(m_1)\rk_1(m_2):.v_i\ot \cdots \ot v_k\in W,
\end{align}
for $m_1,m_2\in \Z, i\in I,$ where the normal ordered product $:\rd_1(m_1)\rk_1(m_2):$ is defined as in \eqref{norpro}.
By using  \eqref{barl} and \eqref{L24}, it is easy to check that
\begin{equation}\begin{split}\label{bardm1}
&L(m).v_i-\bar{\rd}^M(m,1).v_i\\
=\,&L(m).v_i-\bar{L}^M(m,0).v_i+\bar{\rd}_1^M(m,1).v_i\\
=\,&\sum_{m_1,m_2\le M, m_1+m_2=m,m_2\ne 0}\frac{1}{c_i}(1-\frac{m_1}{m_2}):\rd_1(m_1)\rk_1(m_2):.v_i,
\end{split}\end{equation}
for $m\le M$ and $i\in I$.
Then, by taking $s=1$ in \eqref{dms}, we obtain from \eqref{bardm1} and \eqref{rkrd} that
\begin{align*}
v_1\ot \cdots \ot L(m).v_i\ot \cdots \ot v_k\in W,
\end{align*}
for $m\in \Z, i\in I,$ as required.
\end{proof}

\begin{lemt}\label{dhwv} If $\wt\fg(\mu)_+.\bm{w}=0$ for some $\bm{w}\in \wh V(\bm{\lambda},\bm{a})$, then $\bm{w}\in \C v_{\bm{\lambda}}$.
\end{lemt}
\begin{proof}
Let $\{w^{1}_{s}\}_{s\in I_1}$ (respectively $\{w^{2}_{s}\}_{s\in I_2}$, $\cdots$, $\{w^{k}_{s}\}_{s\in I_k}$)
be a basis of $V_{\wh\fg}(\lambda_1)$ (respectively $V_{\wh\fg}(\lambda_2)$, $\cdots$, $V_{\wh\fg}(\lambda_k)$) consisting of weight vectors,
where $I_j$ are index sets for $1\le j\le k$. We assume that $\bm{w}\ne 0$, and then there is a linearly independent set of vectors in $\wh V(\bm{\lambda},\bm{a})$
\[ \{w^1_{s_{1,p}}\otimes w^2_{s_{2,p}}\otimes\cdots\otimes w^k_{s_{k,p}}\mid p=1,\cdots,t, s_{i,p}\in I_i\}\]
such that
\[\bm{w}=\sum_{1\le p\le t}
b_p\, w^1_{s_{1,p}}\otimes w^2_{s_{2,p}}\otimes\cdots\otimes w^k_{s_{k,p}},\]
for some $b_p\in \C^\times.$ We claim that
\begin{align}\label{claim}\sum_{1\le p\le t}
b_p\, w^1_{s_{1,p}}\ot \cdots \ot g.w^i_{s_{i,p}}\ot \cdots \ot w^k_{s_{k,p}}=0,
\end{align}
for all $ g\in \bar\ff_+$, and which will be proved later on. We suppose that this claim holds. Then, for each $i\in I$ and $p=1,2,\cdots,t$, we deduce from \eqref{claim} that
\begin{align}\label{iip}
\sum_{q\in I(i,p)} b_q\, (g.w^i_{s_{i,q}})=g.(\sum_{q\in I(i,p)}b_q\, w^i_{s_{i,q}})=0,
\end{align}
for all $ g\in \bar\ff_+,$
where
\[I(i,p)=\{1\le q\le t\mid s_{j,q}=s_{j,p},\ \forall\, j\in I\setminus\{i\}\}.\]
It is obvious that, if $\bar\ff_+.u_i=0$ for some $u_i\in V_{\bar\ff}(\bar\lambda_i)$, then
  $u_i\in \C v_{\bar\lambda_i}$.
Therefore, from \eqref{iip}, we obtain, for each $i\in I$ and $p=1,2,\cdots,t$, that
\[\sum_{q\in I(i,p)}b_qw^i_{s_{i,q}}\in \C v_{\bar\lambda_i}.\]
This implies that $I(i,p)=\{p\}$ and $w^i_{s_{i,p}}=b_{i,p} v_{\bar\lambda_i}$ for some $b_{i,p}\in \C^\times$.
Thus we have that $t=1$ and $w^i_{s_{i,1}}\in \C v_{\bar\lambda_i}$ for all $i\in I$, as required.

Now we turn to prove the claim \eqref{claim}. We
fix a sufficiently large integer $M$ so that
\[w_{s_{i,p}}^i\in \bigoplus_{n\le M}V_{\bar\ff}(\bar\lambda_i)_n\]
for all  $i\in I$ and $1\le p\le t$. This implies that for all $m>M$, $i\in I$ and $p=1,\cdots,t$,
\[\phi_n^+(m).w^i_{s_{i,p}}=\rk_1(m).w^i_{s_{i,p}}=\rd_1(m).w^i_{s_{i,p}}=\bar{L}(m).w^i_{s_{i,p}}=x(m).w^i_{s_{i,p}}=0,\]
where $n\in \Z$ and $x\in \dg$.
By a similar argument as we did in \eqref{quo1},  we have the following system of equations
\begin{align*}
t_0^mt_1^n\rk_0.\bm{w}=\sum_{i=1}^k\sum_{s=0}^{2M-m} a_i^n n^s(\sum_{1\le p\le t}
b_p\, w^1_{s_{1,p}}\ot \cdots \ot \phi^M(m,s).w^i_{s_{i,p}}\ot \cdots \ot w^k_{s_{k,p}})
\end{align*}
for $m\le M$, and $1\le n\le k(2M-m+1)$. Note that, if $m>0$, then $t_0^mt_1^n \rk_0\in \wt\fg(\mu)_+$ and hence $t_0^mt_1^n\rk_0.\bm{w}=0$.
By solving the above system of linear equations, one gets that
\[\sum_{1\le p\le t}
b_p\, w^1_{s_{1,p}}\ot \cdots \ot \phi^M(m,s).w^i_{s_{i,p}}\ot \cdots \ot w^k_{s_{k,p}}=0,\]
for $ m>0,\, 0\le s\le 2M-m.$
Thus, by taking $s=0$ and using \eqref{L24}, we have that
\begin{align}\label{k1van}
\sum_{1\le p\le t}b_p\, w^1_{s_{1,p}}\ot \cdots \ot \rk_1(m).w^i_{s_{i,p}}\ot \cdots \ot w^k_{s_{k,p}}=0,
\end{align}
for $m>0.$
Furthermore, as we did in the proof of \eqref{quo2}, for every $t_0^m\ot x\in \fg_+$ with $m\le M$, one has that
\begin{align*}
0=t_0^mt_1^n\ot x.\bm{w}=\sum_{i=1}^k\sum_{s=0}^{3M-m} a_i^n n^s(\sum_{1\le p\le t}
b_p\, w^1_{s_{1,p}}\ot \cdots \ot x^M(m,s).w^i_{s_{i,p}}\ot \cdots \ot w^k_{s_{k,p}})
\end{align*}
where $1\le n\le k(3M-m+1)$.
By taking $s=1$, the above identity implies that
\begin{align}\label{xvan}
\sum_{1\le p\le t}b_p\, w^1_{s_{1,p}}\ot \cdots \ot x(m).w^i_{s_{i,p}}\ot \cdots \ot w^k_{s_{k,p}}=0,
\end{align}
for $t_0^m\ot x\in \fg_+.$
Finally, as we did in the proof of \eqref{quo3}, for every $0<m\le M$, we have the following system of equations
\begin{align*}
0=\tilde{\rd}(m,n).\bm{w}=\sum_{i=1}^k\sum_{s=0}^{3M-m+2} a_i^n n^s(\sum_{1\le p\le t}
b_p\, w^1_{s_{1,p}}\ot \cdots \ot \bar{\rd}^M(m,s).w^i_{s_{i,p}}\ot \cdots \ot w^k_{s_{k,p}}),
\end{align*}
for $1\le n\le k(3M-m+3)$.
Then we solve this system of linear equations to obtain, for all  $0<m\le M$ and $0\le s\le k(2M-m+3)$, that
\begin{align}\label{fin}
\sum_{1\le p\le t}
b_p\, w^1_{s_{1,p}}\ot \cdots \ot \bar{\rd}^M(m,s).w^i_{s_{i,p}}\ot \cdots \ot w^k_{s_{k,p}}=0.\end{align}
By taking $s=0,1$ respectively in \eqref{fin}, we obtain
\begin{align}\label{d1van}
\sum_{1\le p\le t}b_p\, w^1_{s_{1,p}}\ot \cdots \ot \rd_1(m).w^i_{s_{i,p}}\ot \cdots \ot w^k_{s_{k,p}}=0,
\end{align}
and
\begin{align}\label{barlvan}
\sum_{1\le p\le t}b_p\, w^1_{s_{1,p}}\ot \cdots \ot L(m).w^i_{s_{i,p}}\ot \cdots \ot w^k_{s_{k,p}}=0,
\end{align}
for $m>0.$
The claim \eqref{claim} then follows from \eqref{k1van}, \eqref{xvan}, \eqref{d1van} and \eqref{barlvan}.
This completes the proof of this lemma.
\end{proof}

Now we are in a position to complete the proof of Theorem \ref{irr}  by applying the previous four lemmas. Recall that $W$ is a non-zero submodule of $\wh V(\bm{\lambda},\bm{a})$ and we need to show that $W=\wh V(\bm{\lambda},\bm{a})$.  Recall from Lemma \ref{basic}(i) that every weight in $W$ has the
form $\underline{\lambda}-\eta$ for some $\eta\in \Gamma_+$. We choose a weight $\underline{\lambda}-\eta_0$ of $W$ which is maximal in the sense that there does not exist any non-zero weight vector in  $W$ with weight $\underline{\lambda}-\eta'$ so that  $\eta_0-\eta'\in \Gamma_+$.
Let $w\in W$ be a non-zero weight vector with weight $\underline{\lambda}-\eta_0$. Then we have $\wt\fg(\mu)_+.w=0$, and hence
$w\in \C v_{\bm{\lambda}}$ by applying Proposition \ref{dhwv}.
Therefore, we obtain that $v_{\bm{\lambda}}\in W$. Moreover, it follows from Lemma \ref{lem1}, Lemma \ref{lem2} and Lemma \ref{lem3} that
\begin{align*}
v_{\bar\lambda_1}\ot \cdots \ot g.v_{\bar\lambda_i}\ot \cdots v_{\bar\lambda_k}\in W,
\end{align*}
for all $g\in \bar\ff$ and $i\in I$. This implies that $W=\wh V(\bm{\lambda},\bm{a})$, and thus we complete the proof of Theorem \ref{irr}.

\end{document}